\documentclass[12pt]{article}
\usepackage{amssymb, amsmath, url, graphicx}

\def\3{\subset }
\def\4{\subseteq }
\def\<{\left<}
\def\>{\right>}

\def\bit{\begin{itemize}}
\def\eit{\end{itemize}}
\def\3{\subset }
\def\4{\subseteq }
\def\ov{\overline}

\def\calm{{\cal M}}
\def\0{\leqno}

\def\barr{\begin{array}}
\def\earr{\end{array}}

\def\Z{{\rlap{$\kern2pt{\rm Z}$}{\rm Z}\,}}

\title{\bf Non-CLT groups of order $pq^3$}
\author{Marius T\u arn\u auceanu}
\date{February 17, 2015}

\begin{document}

\maketitle

\begin{abstract}
In this note we give a characterization of finite groups of order
$pq^3$ ($p$, $q$ primes) that fail to satisfy the Converse of
Lagrange's Theorem.
\end{abstract}

\noindent{\bf MSC (2010):} Primary 20D99;  Secondary 20E99.

\noindent{\bf Key words:} CLT-group, Lagrange's theorem.

\section{Introduction}

All groups considered in this note are finite. A group is said to
be CLT if it possesses subgroups of every possible order (that is,
it satisfies the Converse of Lagrange's Theorem) and non-CLT
otherwise. It is well-known that CLT groups are solvable (see
\cite{7}) and that supersolvable groups are CLT (see \cite{6}).
Recall also that the inclusion between the classes of CLT groups
and solvable groups, as well as the inclusion between the classes
of supersolvable groups and CLT groups are proper (see, for
example, \cite{3}).
\bigskip

An important class of solvable groups, which are not necessarily
supersolvable, consists of the groups of order
$p^{\alpha}q^{\beta}$ with $p, q$ primes and $\alpha, \beta \in
\mathbb{N}$. So, a natural question is whether such a group is
CLT. The Baskaran's papers \cite{1} and \cite{2} answer this
question for the particular cases $\alpha=1,\, \beta=2$ and
$\alpha=\beta=2$, respectively. In the current paper we study the
case $\alpha=1,\, \beta=3$. Our main theorem gives necessary and
sufficient conditions to exist non-CLT groups of order $pq^3$ and
describes the structure of these groups.

\bigskip\noindent{\bf Theorem 1.1.} {\it Let $p$ and $q$ be two primes.
Then there exists a non-CLT group of order $pq^3$ if and only if
either $p$ divides $q+1$ or $p$ divides $q^2+q+1$. Moreover,
excepting the case $p=3,\, q=2$ in which one obtains a unique
non-CLT group, namely $SL(2,3)$, all non-CLT groups of order
$pq^3$ are nontrivial semidirect products of a normal subgroup
$H\cong \mathbb{Z}_q^3$ or $H\cong E(q^3)$ by a subgroup $K\cong
\mathbb{Z}_p$.}
\bigskip

Most of our notation is standard and will usually not be repeated
here. For elementary concepts and results on group theory we refer
to \cite{4} and \cite{5}.

\section{Proofs of the main results}

\noindent{\bf Proof of Theorem 1.1.} Let $G$ be a non-CLT group of
order $pq^3$. Since supersolvable groups are CLT, we easily infer
that:
\begin{itemize}
\item[--] $p\neq q$;
\item[--] $G$ has no normal Sylow $p$-subgroup;
\item[--] $q\not\equiv 1 \hspace{1mm}({\rm mod}\hspace{1mm} p)$.
\end{itemize}Then the number $n_p$ of Sylow $p$-subgroups of $G$
must be $q^2$ or $q^3$. On the other hand, Theorem 1.32 of
\cite{5} shows that $G$ possesses a normal Sylow $q$-subgroup,
except when $\mid G\mid\hspace{1mm}=24$.

\bigskip\noindent{\bf \hspace{10mm}Case 1.} $\mid G\mid\hspace{1mm}=24$

\noindent In this case, by investigating the 15 types of groups of
order 24, we deduce that the only possibility is $G\cong SL(2,3)$
(this is non-CLT because it has no subgroup of order 12).

\bigskip\noindent{\bf \hspace{10mm}Case 2.} $\mid G\mid\hspace{1mm}\neq 24$

\noindent In this case $G$ is a nontrivial semidirect product of a
normal subgroup $H$ of order $q^3$ by a subgroup $K$ of order $p$.
Since the classification of groups of order $q^3$ depends on the
parity of $q$, we distinguish the following two subcases.

\bigskip\noindent{\bf \hspace{20mm}Subcase 2.1.} $q=2$

\noindent The conditions $n_p\in\{4,8\}$ and $n_p\equiv 1
\hspace{1mm}({\rm mod}\hspace{1mm} p)$ lead to $p=7$, $n_p=8$ and
$\mid G\mid\hspace{1mm}=56$. Up to isomorphism, there is a unique
group of order 56 without normal Sylow 7-subgroups, and the Sylow
2-subgroup of this group is elementary abelian. Moreover, we can
easily check that it does not possess subgroups of order 28, i.e.
it is indeed non-CLT.

\bigskip\noindent{\bf \hspace{20mm}Subcase 2.2.} $q\neq 2$

\noindent The groups of order $q^3$ for $q$ odd are either
abelian, namely $\mathbb{Z}_{q^3}$,
$\mathbb{Z}_q\times\mathbb{Z}_{q^2}$ and $\mathbb{Z}_q^3$, or
nonabelian, namely $M(q^3)=\langle x,y \mid x^{q^2}=y^q=1,
y^{-1}xy=x^{q+1}\rangle$ and $E(q^3)=\langle x,y \mid
x^q=y^q=[x,y]^q=1, [x,y]\in Z(E(q^3))\rangle$. We observe that we
cannot have $H\cong\mathbb{Z}_{q^3}$, because in this case $G$
would be metacyclic and therefore CLT. The number of automorphisms
of the other four groups is:
\begin{itemize}
\item[--] $\mid {\rm Aut}(\mathbb{Z}_q\times\mathbb{Z}_{q^2})\mid\hspace{1mm}=q^3(q-1)^2$,
\item[--] $\mid {\rm Aut}(\mathbb{Z}_q^3)\mid\hspace{1mm}=q^3(q-1)(q^2-1)(q^3-1)$,
\item[--] $\mid {\rm Aut}(M(q^3))\mid\hspace{1mm}=q^3(q-1)^2$,
\item[--] $\mid {\rm Aut}(E(q^3))\mid\hspace{1mm}=q^3(q-1)^2(q+1)$.
\end{itemize}Since there is a nontrivial homomorphism from $K\cong\mathbb{Z}_p$ to ${\rm Aut}(H)$, $p$ must
divide $\mid {\rm Aut}(H)\mid$, which implies that either $H\cong
\mathbb{Z}_q^3$ or $H\cong E(q^3)$. It is also clear that one of
the conditions $p\mid q+1$ or $p\mid q^2+q+1$ is verified.
\bigskip

Conversely, suppose first that $p\mid q^2+q+1$. Then every
nontrivial semidirect product $G$ of a normal elementary abelian
subgroup of order $q^3$ by a subgroup of order $p$ is non-CLT.
Indeed, if we assume that $G$ possesses a subgroup of order
$pq^2$, say $G_1$, then there is a Sylow $p$-subgroup $S_p$ of $G$
such that $S_p\subset G_1$. By applying the Sylow's theorems for
$G_1$, it follows that $S_p$ is normal in $G_1$, that is
$G_1\subseteq N_G(S_p)$. This leads to $N_G(S_p)=G_1$ and
therefore $q^3=n_p=[G:N_G(S_p)]=[G:G_1]=q$, a contradiction.

Suppose next that $p\mid q+1$ and let $G$ be a nontrivial
semidirect product of a normal subgroup isomorphic to $E(q^3)$ by
a subgroup $\langle a\rangle$ of order $p$ such that $a$ commutes
with $[x,y]$ ($x$ and $y$ denote the generators of $E(q^3)$, as
above). We will prove that $G$ is non-CLT by showing again that it
does not possess subgroups of order $pq^2$. If $G_1$ is such a
subgroup and $S_p$ is a Sylow $p$-subgroup of $G$ contained in
$G_1$, then we can assume that $S_p=\langle a\rangle$. On the
other hand, it is obvious that $\langle [x,y]\rangle=\Phi(\langle
x,y\rangle)\subset G_1$. Then $G_1$ contains the commuting
subgroups $\langle a\rangle\cong\mathbb{Z}_p$ and $\langle
[x,y]\rangle\cong\mathbb{Z}_q$. Consequently, it has subgroups of
order $pq$, i.e. it is CLT. By the main theorem of \cite{1}, we
deduce that $G_1$ is necessarily abelian, which implies
$G_1\subseteq N_G(S_p)$. As in the first part of this implication,
one obtains $N_G(S_p)=G_1$ and hence
$q^2=n_p=[G:N_G(S_p)]=[G:G_1]=q$, a contradiction. This completes
the proof.
\hfill\rule{1,5mm}{1,5mm}
\bigskip

An immediate consequence of Theorem 1.1 is given by the following
corollary.

\bigskip\noindent{\bf Corollary 2.1.} {\it The unique non-CLT groups of order $8p$
with $p$ prime are $SL(2,3)$ and the group of order {\rm 56}
described above.}
\bigskip

Next, let $G$ be a non-CLT group of order $pq^3$ ($p$, $q$
primes), $S_q$ be a Sylow $q$-subgroup of $G$ and $\calm$ be the
set of subgroups of order $q^2$ of $S_q$. We consider the
conjugation action of $G$ on $\calm$ and we choose a set of
representatives $\{H_1,H_2,...,H_k\}$ for the conjugacy classes.
Clearly, we have $S_q\subseteq N_G(H_i)$,
$\forall\hspace{1mm}i=\ov{1,k}$. On the other hand, every $H_i$ is
not normal in $G$ because $G$ does not possess subgroups of order
$pq^2$. These prove that $N_G(H_i)=S_q$ and therefore
$[G:N_G(H_i)]=p$, $\forall\hspace{1mm}i=\ov{1,k}$. Thus, we infer
that $p$ divides $\mid\calm\mid$.
\bigskip

Since the numbers of subgroups of order $q^2$ of $E(q^3)$ and
$\mathbb{Z}_q^3$ are $q+1$ and $q^2+q+1$, respectively, the above
remark shows that Theorem 1.1 can be reformulated in the following
way for $q\neq 2$.

\bigskip\noindent{\bf Corollary 2.2.} {\it Given two primes $p$ and $q$, with $q\neq 2$, the following
statements are true:
\begin{itemize}
\item[{\rm a)}] There exists a non-CLT group of order $pq^3$ having a Sylow $q$-subgroup isomorphic to $E(q^3)$ if and only if $p$ divides $q+1$.
\item[{\rm b)}] There exists a non-CLT group of order $pq^3$ having an elementary abelian Sylow $q$-subgroup if and only if $p$ divides $q^2+q+1$.
\end{itemize}}
\bigskip

Finally, we remark that a result similar with Corollary 2.2 also
holds in the case $q=2$ (notice that the Sylow 2-subgroups of
$SL(2,3)$ are isomorphic to the well-known quaternion group
$Q_8$).

\bigskip\noindent{\bf Corollary 2.3.}  {\it Given a prime $p$, the following
statements are true:
\begin{itemize}
\item[{\rm a)}] There exists a non-CLT group of order $8p$ having a Sylow {\rm 2}-subgroup isomorphic to $Q_8$ if and only if $p=3$.
\item[{\rm b)}] There exists a non-CLT group of order $8p$ having an elementary abelian Sylow {\rm 2}-subgroup if and only if $p=7$.
\end{itemize}}

\vspace*{5ex}\small

\hfill
\begin{minipage}[t]{5cm}
Marius T\u arn\u auceanu \\
Faculty of  Mathematics \\
``Al.I. Cuza'' University \\
Ia\c si, Romania \\
e-mail: {\tt tarnauc@uaic.ro}
\end{minipage}


\begin{thebibliography}{10}
\bibitem{1} Baskaran, S., {\it CLT and non-CLT groups}, I, Indian J. Math. {\bf 14} (1972), 81-82.
\bibitem{2} Baskaran, S., {\it CLT and non-CLT groups of order $p^2q^2$}, Fund. Math. {\bf 92} (1976), 1-7.
\bibitem{3} Bray, H.G., {\it A note on CLT groups}, Pacific J. Math. {\bf 27} (1968), 229-231.
\bibitem{4} Huppert, B., {\it Endliche Gruppen}, I, Springer Verlag, Berlin-Heidelberg-New York, 1967.
\bibitem{5} Isaacs, I.M., {\it Finite group theory}, Amer. Math. Soc., Providence, R.I., 2008.
\bibitem{6} McCarthy, D.J., {\it A survey of partial coverses to Lagrange's theorem on finite groups}, Trans. N.Y. Acad. Sci. {\bf 33} (1971), 586-594.
\bibitem{7} McLain, D.H., {\it The existence of subgroups of given order in finite groups}, Proc. Cambridge Philos. Soc. {\bf 53} (1957), 278-285.
\end{thebibliography}
\end{document}